\documentclass[MSNbibl,number,citesort,seceqn,dvips]{arxbj}
\usepackage{graphicx}
\aid{0}
\volume{20}
\issue{4}
\pubyear{2014}
\firstpage{1717}
\lastpage{1737}
\doi{10.3150/13-BEJ538} %kopijuoti is PTS

\makeatletter

\newcommand{\RMe}{\mathrm{e}}

\newcommand{\mrmd}{\,\mathrm{d}}

\newcommand{\cal}{\mathcal}

\newtheorem{prop}{Proposition}[section]
\newtheorem{theo}{Theorem}[section]
\newtheorem{cor}{Corollary}[section]
\newremark{rem}{Remark}[section]
\newproclaim{defi}{Definition}[section]
\newremark{ex}{Example}[section]

\newcommand{\R}{\mathbb{R}}
\gdef\upqed{\vskip-\belowdisplayskip\vskip-\baselineskip\kern0pt}

\makeatother

\begin{document}
\begin{frontmatter}

\title{The generalized Pareto process; with a view towards application
and simulation}
\runtitle{The generalized Pareto process}

\begin{aug}
%%%% inicialai - be tarpu
\author[1,3]{\inits{A.}\fnms{Ana} \snm{Ferreira}\corref{}\thanksref{1,3}\ead[label=e1]{anafh@isa.utl.pt}} \and
\author[2,3]{\inits{L.}\fnms{Laurens} \snm{de Haan}\thanksref{2,3}\ead[label=e2]{ldehaan@ese.eur.nl}}
%%\runauthor{} %% auto
\address[1]{ISA, Univ Tecn Lisboa, Tapada da Ajuda 1349-017 Lisboa,
Portugal. \printead{e1}}
\address[2]{Erasmus University Rotterdam, P.O. Box 1738, 3000 DR
Rotterdam, The Netherlands.\\ \printead{e2}}
\address[3]{CEAUL, FCUL, Bloco C6 -- Piso 4 Campo Grande, 749-016
Lisboa, Portugal}
\end{aug}

% HISTORY:
\received{\smonth{11} \syear{2012}}
\revised{\smonth{4} \syear{2013}}

% ABSTRACT
%
\begin{abstract}
In extreme value statistics, the peaks-over-threshold method is widely
used. The method is based on the generalized Pareto distribution
characterizing probabilities of exceedances over high thresholds in
$\R^d$. We present a generalization of this concept in the space of
continuous functions. We call this the generalized Pareto process.
Differently from earlier papers, our definition is not based on a
distribution function but on functional properties, and does not need a
reference to a related max-stable process.

As an application, we use the theory to simulate wind fields connected
to disastrous storms on the basis of observed extreme but not
disastrous storms. We also establish the peaks-over-threshold approach
in function space.
\end{abstract}

% KEYWORDS
% visi is mazosios raides ir pagal abecele
%
\begin{keyword}
\kwd{domain of attraction}
\kwd{extreme value theory}
\kwd{functional regular variation}
\kwd{generalized Pareto process}
\kwd{max-stable processes}
\kwd{peaks-over-threshold}
\end{keyword}

\end{frontmatter}

%s1 #&#
\section{Introduction}\label{Introdsect}
Let $C(S)$ be the space of continuous real functions on $S$, equipped
with the supremum norm, where $S$ is a compact subset of $\R^d$. A
stochastic process $X$ in $C(S)$ is in the domain of attraction of some
max-stable process $Y$ if there are continuous functions $a_n(s)$
positive and $b_n(s)$ on $S$ such that the processes
\[
\biggl\{\max_{1\leq i\leq n}\frac{X_i(s)-b_n(s)}{a_n(s)} \biggr\}
_{s\in S}
\]
with $X,X_1,\ldots,X_n$ independent and identically distributed
(i.i.d.), converge in distribution to $Y$ in $C(S)$. Necessary and
sufficient conditions are: uniform convergence of the marginal
distributions and a convergence of measures (in fact a form of regular
variation):
%
%e1.1 #&#
%
\begin{equation}
\label{nudefintrod} \lim_{t\to\infty} tP
\bigl(T_tX(\cdot)\in A \bigr)=\nu(A),
\end{equation}
where:
\begin{itemize}
\item$T_tX(s)= (1+\gamma(s)\frac{X(s)-b_t(s)}{a_t(s)}
)_+^{1/\gamma(s)}$ for all $s\in S$ (with the notation $x_+=\max
(0,x)$ for any real $x$),
\item$\nu$ is a homogeneous measure of order $-1$ on $C^+(S):=\{f\in
C(S)\dvt f\geq0\}$ and,
\item$A$ is any Borel subset of $C^+(S)$ satisfying $\inf\{\sup_{s\in
S} f(s)\dvt f\in A\}>0$ and $\nu(\partial A )=0$.
\end{itemize}
Cf. de Haan and Lin \cite{HaanLin} and de Haan and Ferreira \cite
{HaanFerreira}, Section 9.5. Although in these references only the case
$C=[0,1]$ has been worked out, all arguments are valid for any compact
subset of~$\R^d$ (as remarked in these references). The functions
$a_t(s)$ and $b_t(s)$ are chosen in such a way that the marginal
distributions are in standard form,
\[
\lim_{t\to\infty} tP \biggl( \frac{X(s)-b_t(s)}{a_t(s)}>x \biggr)= \bigl(1+
\gamma(s)x \bigr)^{-1/\gamma(s)},\qquad 1+\gamma(s)x>0,
\]
for all $x\in\R$ and $s\in S$. Here $\gamma$ is a continuous
function. In particular, one may take $b_t(s):=\inf\{x\dvt P(X(s)\leq
x)\geq1-1/t\}$. This is how we choose $b_t(s)$ from now on. One
possible choice of $a_t(s)$ is $a_t(s):=\gamma
(s)(b_{2t}(s)-b_t(s))/(2^{\gamma(s)}-1)$.

From (\ref{nudefintrod}), it follows that
\[
\frac{P ( (1+\gamma(\cdot)({X(\cdot)-b_t(\cdot
)})/{a_t(\cdot)} )_+^{1/\gamma(\cdot)}\in A )}{P
(\sup_{s\in S} (({X(s)-b_t(s)})/{a_t(s)})>0 )}
\]
converges, as $t\to\infty$, and so does
\[
P \biggl( \biggl(1+\gamma(\cdot)\frac{X(\cdot)-b_t(\cdot)}{a_t(\cdot
)} \biggr)_+^{1/\gamma(\cdot)}\in A \Big|
\sup_{s\in S} \frac
{X(s)-b_t(s)}{a_t(s)}>0 \biggr).
\]
The limit constitutes a probability distribution on $C^+(S)$.

This reasoning is similar to how one obtains the generalized Pareto
distributions in $\R$ (Pickands~\cite{Pickands75}; Balkema and
de Haan \cite{BalkemadeHaan}) and in $\R^d$ (Rootz{\'e}n and Tajvidi
\cite{RootzenTajvidi}; Falk, H{\"u}sler and Reiss \cite
{FalkHuslerReiss}). It leads to what we call generalized Pareto processes.

The paper is organized as follows.
The Pareto processes will be dealt with in Section \ref{SPpsect}. As
in the finite-dimensional context, it is convenient to study first
generalized Pareto processes in a standardized form. This is done in
Section \ref{SSPpsect}. The general process is discussed in
Section~\ref{GSPpsect}. In Section \ref{MSPpsect}, a discrete
version of our approach is discussed leading to simple multivariate
Pareto random vectors. The domain of attraction is discussed in
Section \ref{domattrsect}.
In Section \ref{applicsect}, we show that by using the stability
property of generalized Pareto processes one can create extreme storm
fields starting from independent and identically observations of storm
fields. We also establish the peaks-over-threshold approach in function space.

In the following, operations like $w_1+w_2$ or $w_1\wedge w_2$ with
$w_1,w_2\in C(S)$ mean, respectively, $\{w_1(s)+w_2(s)\}_{s\in S}$ and
$\{
w_1(s)\wedge w_2(s)\}_{s\in S}$. Then, with abuse of notation,
operations like $w+x$ or $w\wedge x$ with $w\in C(S)$ and $x\in\R$
mean, respectively, $\{w(s)+x\}_{s\in S}$ and $\{w(s)\wedge x\}_{s\in
S}$. Similarly for products and powers. Then, for example, we shall
simply write $ (1+\gamma\frac{X-b_t}{a_t} )^{1/\gamma}$
for $ \{ (1+\gamma(s)\frac{X(s)-b_t(s)}{a_t(s)}
)^{1/\gamma(s)} \}_{s\in S}$, with $X=\{X(s)\}_{s\in S}$, $a_t=\{
a_t(s)\}_{s\in S}$, $b_t=\{b_t(s)\}_{s\in S}$ and $\gamma=\{\gamma
(s)\}_{s\in S}$.

Denote the class of Borel subsets of a metric space by ${\cal B}
(\cdot)$.

%s2 #&#
\section{Pareto processes}
\label{SPpsect}

%s2.1 #&#
\subsection{The simple Pareto process}
\label{SSPpsect}

Again, let $C^+(S)$ be the space of non-negative real continuous
functions on $S$, with $S$ some compact subset of $\R^d$.

%
%th2.1 #&#
%
\begin{theo}\label{Wtheor}
Let $W$ be a stochastic process in $C^+(S)$ and $\omega_0$ a positive
constant. The following three statements are equivalent:
\begin{enumerate}[3.]
\item(Peaks-over-threshold):
\begin{enumerate}[(a)]
\item[(a)] The expectation $E (W(s)/\sup_{u\in S} W(u) )$ is
positive for all $s\in S$,\vspace*{1pt}
\item[(b)]$P (\sup_{s\in S}W(s)/\omega_0>x )=x^{-1}$, for
$x>1$ (standard Pareto distribution),
\item[(c)]
%
%e2.1 #&#
%
\begin{equation}
\label{Qdeff} P \biggl(\frac{\omega_0 W}{\sup_{s\in S} W(s)}\in B \Big| \sup
_{s\in S}
W(s)>r \biggr)=P \biggl(\frac{\omega_0 W}{\sup_{s\in S}
W(s)}\in B \biggr)
\end{equation}
for all $r>\omega_0$ and $B\in{\cal B} (\bar C_{\omega
_0}^{+}(S) )$ with
%
%e2.2 #&#
%
\begin{equation}
\label{barQspace} \bar C_{\omega_0}^{+}(S):=\Bigl\{f \in C^+(S)\dvt
\sup_{s\in S} f(s)=\omega_0\Bigr\}.
\end{equation}
\end{enumerate}

\item(Random functions):
\begin{enumerate}[(a)]
\item[(a)]$P (\sup_{s\in S} W(s)> \omega_0 )=1$,
\item[(b)]$E (W(s)/\sup_{u\in S} W(u) )>0$ for all $s\in S$,
\item[(c)]
%
%e2.3 #&#
%
\begin{equation}
\label{GPhomogeneity} P(W\in rA)=r^{-1}P(W\in A)
\end{equation}
for all $r>1$ and $A\in{\cal B} (C_{\omega_0}^{+}(S) )$,
where $rA$ means the set $\{rf, f\in A\}$, and
%
%e2.4 #&#
%
\begin{equation}
\label{Qspace} C_{\omega_0}^{+}(S):=\Bigl\{f \in C^+(S)\dvt\sup
_{s\in S} f(s)\geq\omega_0\Bigr\}.
\end{equation}
\end{enumerate}

\item(Constructive approach) $W(s)=YV(s)$, for all $s\in S$, for some
$Y$ and $V=\{V(s)\}_{s\in S}$ satisfying:
\begin{enumerate}[(a)]
\item[(a)]$V\in C^+(S)$ is a stochastic process satisfying $\sup_{s\in S}
V(s)=\omega_0$ a.s., and $EV(s)>0$ for all $s\in S$,
\item[(b)]$Y$ is a standard Pareto random variable, $P(Y\leq y)=1-1/y$, $y>1$,
\item[(c)]$Y$ and $V$ are independent.
\end{enumerate}
\end{enumerate}
\end{theo}

%
%de2.1 #&#
%
\begin{defi}\label{Wdefi}
The process $W$ characterized in Theorem \ref{Wtheor}, with threshold
parameter $\omega_0$, is called simple Pareto process. The probability
measure in (\ref{Qdeff}), that is,
%
%e2.5 #&#
%
\begin{equation}
\label{Qdef} \rho(B)=P \biggl(\frac{\omega_0 W}{\sup_{s\in S} W(s)}\in
B \biggr)\qquad \mbox{for } B\in{\cal B} \bigl(\bar C_{\omega
_0}^{+}(S) \bigr)
\end{equation}
is called the spectral measure.
\end{defi}

Some easy consequences of Theorem \ref{Wtheor}(3) are the
following. The process $W$ is stationary if and only if $V$ is
stationary. Independence at any two points $s_1,s_2\in S$, that is,
$W(s_1)$ and $W(s_2)$ being independent, is not possible. Complete
dependence is equivalent to $V\equiv\omega_0$ a.s. We shall come back
to some of these issues.

\begin{pf*}{Proof of Theorem \ref{Wtheor}}
We start by proving that 1 implies 3. By compactness
and continuity, $\sup_{s\in S}W(s)<\infty$ a.s. Take:
\[
Y=\frac{\sup_{s\in S} W(s)}{\omega_0} \quad\mbox{and}\quad V=\frac{\omega_0
W}{\sup_{s\in S} W(s)}.
\]
Then (a), (b) and (c) are straightforward.

Next, we prove that 3 implies 2. Let
\[
A_{r,B}= \biggl\{f \in C^+(S)\dvt\sup_{s\in S} f(s)/
\omega_0>r, \frac
{\omega_0 f}{\sup_{s\in S} f(s)}\in B \biggr\}= r A_{1,B}
\]
for all $r>1$ and $B\in{\cal B} (\bar C_{\omega_0}^{+}(S)
)$. Then,
\begin{eqnarray*}
P (W\in A_{r,B} )&=&P \biggl(\sup_{s\in S} W(s)/
\omega_0>r, \frac{\omega_0 W}{\sup_{s\in S} W(s)}\in B \biggr)
\\
&=&P(Y>r, V\in B)=P(Y>r) P(V\in B)
\\
&=&\frac1 r P \biggl(\sup_{s\in S} W(s)/\omega_0>1,
\frac{\omega_0
W}{\sup_{s\in S} W(s)}\in B \biggr) \\
&=&\frac1 r P (W\in A_{1,B} )
\end{eqnarray*}
using in particular the independence of $Y$ and $V$ and $P (\sup_{s\in
S} W(s)/\omega_0>1 )=1$. Since $P(rA)=r^{-1}P(A)$ holds
for any of the above sets, it must also hold for all Borel sets in the
statement.

Finally, check that 2 implies 1. For any $r> 1$, by
(c) and (a),
\[
P \biggl(\frac{\sup_{s\in S} W(s)}{\omega_0}> r \biggr)=\frac1 r P
\biggl(\frac{\sup_{s\in S} W(s)}{\omega_0}> 1
\biggr)=\frac1 r.
\]
Also for any $B\in{\cal B} (\bar C_{\omega_0}^{+}(S) )$,
\begin{eqnarray*}
P \biggl(\sup_{s\in S} W(s)/\omega_0> r,
\frac{\omega_0 W}{\sup_{s\in S} W(s)}\in B \biggr)
&=&\frac1 r P \biggl(\sup_{s\in S} W(s)/\omega_0> 1,
\frac{\omega_0
W}{\sup_{s\in S} W(s)}\in B \biggr) \\
&=&\frac1 r P \biggl(\frac{\omega_0
W}{\sup_{s\in S} W(s)}\in B \biggr)
\end{eqnarray*}
since $\sup_{s\in S} W(s)>\omega_0$ holds a.s. That is, it follows
that $\sup_{s\in S} W(s)/\omega_0$ is univariate Pareto distributed
and, $\sup_{s\in S} W(s)$ and $W/\sup_{s\in S} W(s)$ are independent.
\end{pf*}

The following properties are direct consequences:

%
%co2.1 #&#
%
\begin{cor}\label{supxcor}
For any simple Pareto process $W$, the random variable
$\omega_0^{-1} \sup_{s\in S} W(s)$ has standard Pareto distribution.
\end{cor}

%
%co2.2 #&#
%
\begin{cor} $W\in C^+(S)$ is a simple Pareto process if and only if any
of the two equivalent statements hold:
\begin{enumerate}
\item
\begin{enumerate}[(a)]
\item[(a)]$E (W(s)/\sup_{u\in S} W(u) )>0$ for all $s\in S$,
\item[(b)]$P (\sup_{s\in S}W(s)/\omega_0>x )=x^{-1}$, for $x>1$,
\item[(c)]
%
%e2.6 #&#
%
\begin{equation}
\label{GPmu} P \Bigl(W\in rA \big|\sup_{s\in S} W(s)>r
\omega_0 \Bigr)=P (W\in A )
\end{equation}
for all $r>1$ and $A\in{\cal B} (C_{\omega_0}^{+}(S) )$.
\end{enumerate}
\item
\begin{enumerate}[(a)]
\item[(a)]$E (W(s)/\sup_{u\in S} W(u) )>0$ for all $s\in S$,
\item[(b)]
%
%e2.7 #&#
%
\begin{equation}
\label{GPrv} P \biggl(\sup_{s\in S}
\frac{W(s)}{\omega_0}>r, \frac{\omega_0
W}{\sup_{s\in S} W(s)}\in B \biggr)=\frac{\rho(B)}r
\end{equation}
for all $r>1$ and $B\in{\cal B} (\bar C_{\omega_0}^{+}(S) )$.
\end{enumerate}
\end{enumerate}
\end{cor}

From (\ref{GPmu}), we see that the probability distribution of $W$
serves in fact as the exponent measure in max-stable processes (cf.
de Haan and Ferreira \cite{HaanFerreira}, Section 9.3).
Characterization 2 suggests ways for testing and modeling
Pareto processes.

We proceed to express the distribution function of $W$ in terms of the
probability distribution of $V$ from Theorem \ref{Wtheor}(3)
and Definition \ref{Wdefi}.

Let $w,W\in C^+(S)$. The notation $W\leq w$ will mean $W(s)\leq w(s)$
for all $s\in S$. Similarly for $W>w$ and $W\nleq w$. Clearly, the
latter two are not the same.

Take for the conditional expectation,
\[
E \bigl(g(V)|V\in B \bigr)=\frac1{\rho(B)}\int_B g(v) \mrmd
\rho(v),\qquad B\in{\cal B} \bigl(\bar C_{\omega_0}^{+}(S) \bigr),
\]
defined in the usual sense and whenever $\rho(B)=P(V\in B)>0$, with
$g$ a real functional (e.g., see Billingsley \cite{Billingsley95},
Section 34).

%
%pr2.1 #&#
%
\begin{prop}\label{Wdf1prop} Let $w,W\in C^+(S)$, with $W$ a simple
Pareto process. Let $S_0=\{s\in S\dvt w(s)=0\}$, $\bar S_0=S\setminus S_0$
the complement of $S_0$, and $B_0=\{f\in\bar C_{\omega_0}^+(S)\dvt\inf
_{s\in\bar S_0} \frac{w(s)}{f(s)}\geq1$
and $f(s)=0$ for $s\in S_0\}$. Then
%
%e2.8 #&#
%
\begin{equation}
\label{Wdf} P(W\leq w)=\cases{ \displaystyle \rho(B_0) \biggl\{1-E \biggl(\sup
_{s\in\bar S_0}\frac{V(s)}{w(s)} \Big| V\in B_0 \biggr) \biggr
\}, &\quad if $\rho(B_0)>0$,
\cr
0, &\quad if $\rho(B_0)=0$.}
\end{equation}
\end{prop}

\begin{pf}
\begin{eqnarray*}
P(W\leq w)&=&P\bigl(W(s)\leq w(s) \mbox{ for } s\in\bar S_0
\mbox{ and } W(s)\leq w(s) \mbox{ for } s\in S_0\bigr)
\\
&=&P \biggl(Y\leq\inf_{s\in\bar S_0} \frac{w(s)}{V(s)} \mbox{ and }
V(s)=0 \mbox{ for } s\in S_0 \biggr)
\\
&=&P \biggl(Y\leq\inf_{s\in\bar S_0} \frac{w(s)}{V(s)} \mbox{ and }
V(s)=0 \mbox{ for } s\in S_0 \mbox{ and } \inf_{s\in\bar S_0}
\frac
{w(s)}{V(s)}\geq1 \biggr)
\\
&&{} + P \biggl(Y\leq\inf_{s\in\bar S_0} \frac{w(s)}{V(s)} \mbox{ and }
V(s)=0 \mbox{ for } s\in S_0 \mbox{ and } \inf_{s\in\bar S_0}
\frac{w(s)}{V(s)}< 1 \biggr)
\\
&=&\int_{B_0} P \biggl(Y\leq\inf_{s\in\bar S_0}
\frac
{w(s)}{v(s)} \biggr)\mrmd \rho(v)
\\
&=&\int_{B_0} 1- \sup_{s\in\bar S_0}
\frac{v(s)}{w(s)}\mrmd \rho(v)=\rho(B_0)-\int_{B_0}
\sup_{s\in\bar S_0} \frac
{v(s)}{w(s)}\mrmd \rho(v)
\\
&=&\rho(B_0) \biggl\{1-E \biggl(\sup_{s\in\bar S_0}
\frac
{V(s)}{w(s)} \Big| V\in B_0 \biggr) \biggr\},
\end{eqnarray*}
where the last but two equality follows by the fact that the second
summand in the previous equality is zero and from the independence of
$Y$ and $V$.
\end{pf}

%
%co2.3 #&#
%
\begin{cor}
Under the conditions of Proposition \ref{Wdf1prop},
%
%e2.9 #&#
%
\begin{equation}
\label{fdrho1} P(W\leq w)=1-E \biggl(\sup_{s\in\bar S_0}
\frac{V(s)}{w(s)} \biggr) \qquad\mbox{if } \rho(B_0)=1.
\end{equation}
\end{cor}

The following is obtained in the particular case of $w$ being strictly positive.

%
%pr2.2 #&#
%
\begin{prop}\label{Wdfprop} Let $w,W\in C^+(S)$, with $w$ positive
and $W$ a simple Pareto process. Then
%
%e2.10 #&#
%
\begin{equation}
\label{Wdf>0} P(W\leq w)=E \biggl(\sup_{s\in S}\frac{V(s)}{w(s)\wedge
\omega
_0}
\biggr)-E \biggl(\sup_{s\in S}\frac{V(s)}{w(s)} \biggr).
\end{equation}
%
%In particular, $P(W\not\leq\omega_0)=1$.
\end{prop}

\begin{pf}%[Proof of Proposition \ref{Wdfprop}]
(i) First, consider the case $\inf_{s\in S} w(s)\geq\omega_0$. Use
Theorem \ref{Wtheor}, part 3,
%
%e2.11 #&#
%
\begin{equation}
\label{GPcartesiandf1} P(W\leq w)=P(YV\leq w)=P \biggl(Y\leq\inf
_{s\in S}\frac
{w(s)}{V(s)} \biggr)=1-E \biggl(\sup
_{s\in S}\frac{V(s)}{w(s)} \biggr)
\end{equation}
hence,
%
%e2.12 #&#
%
\begin{equation}
\label{GPcartesiandf} P(W\nleq w)%=\int_{C_{\omega_0}^{+}(\R)}\sup_{s\in
% dQ(v)
=E \biggl(\sup
_{s\in S}\frac{V(s)}{w(s)} \biggr).%=\omega_0 E\left(
\end{equation}

(ii) The probability measure of $W$ on $C_{\omega_0}^{+}(S)$ can be
extended to a measure $\nu$ on $C^+(S)$ while keeping the homogeneity
relation (\ref{GPhomogeneity}) as follows: for any Borel set $B$ such that
\[
\sup_{f\in B}\sup_{s\in S} f(s) \leq
\omega_0 \quad\mbox{and}\quad 0<\varepsilon<\inf_{f\in B}\sup
_{s\in S} f(s),
\]
we define
\[
\nu(B):=\frac{\omega_0}{\varepsilon}P \biggl(W\in\frac{\omega
_0}{\varepsilon} B \biggr).
\]
This measure (the same as in (\ref{nudefintrod})) is homogeneous of
order $-1$:
\[
\nu(rB)=r^{-1}\nu(B) \qquad\mbox{for all } r>0 \quad\mbox{and}\quad B\in{\cal B}
\bigl(C^+(S) \bigr).
\]
Then, the probability distribution of $W$ is the restriction of $\nu$
to $C_{\omega_0}^{+}(S)$, that is, for $B\in{\cal B}
(C^+(S) )$,
%
%e2.13 #&#
%
\begin{equation}
\label{PWaux} P(W\in B)=\nu\Bigl\{f\in B,\sup_{s\in S} f(s) >
\omega_0 \Bigr\}.
\end{equation}

Hence, by the homogeneity property of $\nu$, (\ref{PWaux}) and (\ref
{GPcartesiandf}) in that order:
%
%e2.14 #&#
%
\begin{eqnarray}
\label{nu++E}
\nu\{f\nleq w \}&=&\frac{\omega_0}{\inf_{s\in S} w(s)}\nu\biggl\{
f\nleq
\frac{w \omega_0}{\inf_{s\in S} w(s)} \biggr\}
\nonumber\\[-8pt]\\[-8pt]
&=&\frac{\omega_0}{\inf_{s\in S} w(s)}P \biggl(W\nleq\frac{w \omega
_0}{\inf_{s\in S} w(s)} \biggr)=E \biggl(\sup
_{s\in S}\frac
{V(s)}{w(s)} \biggr).
\nonumber
\end{eqnarray}
By (\ref{PWaux}), elementary set-measure operations and (\ref{nu++E})
in that order:
\begin{eqnarray*}
P(W\nleq w)&=&\nu\{f\nleq w,f\nleq\omega_0 \}
\\
&=&\nu\{f\nleq w \}+\nu\{f\nleq\omega_0 \} -\nu\{f\nleq w \mbox{
or } f\nleq\omega_0 \}
\\
&=&\nu\{f\nleq w \}+\nu\{f\nleq\omega_0 \} -\nu\{f\nleq w\wedge
\omega_0 \}
\\
&=&E \biggl(\sup_{s\in S}\frac{V(s)}{w(s)} \biggr)+1-E \biggl(
\sup_{s\in S}\frac{V(s)}{w(s)\wedge\omega_0} \biggr).
\end{eqnarray*}
\upqed\end{pf}

Note that $E (\sup_{s\in\bar S_0}\frac{V(s)}{w(s)\wedge\omega
_0} | V\in B_0 )=1$, whenever $\rho(B_0)>0$, which links the
results of Propositions \ref{Wdf1prop} and \ref{Wdfprop}.\vadjust{\goodbreak}

The following formulas might also be useful.

%
%co2.4 #&#
%
\begin{cor}\label{Wdfcor1}
Let $w, W\in C^+(S)$, with $W$ a simple Pareto process. Then:
\begin{enumerate}[(a)]
\item[(a)] With $B_1=\{f\in\bar C_{\omega_0}^+(S)\dvt\sup_{s\in S}
\frac{w(s)}{f(s)}>1$ and $\inf_{s\in S} f(s)>0\}$,
%
%e2.15 #&#
%
\begin{equation}
\label{GPcartesiandf>} P(W> w)= \cases{ \displaystyle \rho(B_1)E \biggl(\inf
_{s\in S}\frac{V(s)}{w(s)} \Big| V\in B_1 \biggr), &\quad if $
\rho(B_1)>0$,
\cr
0, &\quad if $\rho(B_1)=0$.}
\end{equation}
In particular, if $P(V>0)>0$ and $\sup_{s\in S} w(s)>\omega_0$,
%
%e2.16 #&#
%
\begin{equation}
\label{GPcartesiandf>a} P(W> w)=P(V>0)E \biggl(\inf_{s\in S}
\frac{V(s)}{w(s)} \Big| V>0 \biggr).
\end{equation}
\item[(b)] If $w>0$ and $\sup_{s\in S} w(s)>\omega_0$,
%
%e2.17 #&#
%
\begin{equation}
\label{P>w>0} P(W> w)=E \biggl(\inf_{s\in S}\frac{V(s)}{w(s)}
\biggr).
\end{equation}
\item[(c)] If $E (\inf_{s\in S}V(s) )>0$, for $x\in\R$,
%
%e2.18 #&#
%
\begin{equation}
\label{GPcartesiandf>b} P(W> x|W>\omega_0)=\cases{ 1, &\quad $x \leq
\omega_0$,
\cr
\omega_0/x, &\quad $x >\omega_0$.}
\end{equation}
\item[(d)] If $E (\inf_{s\in S}V(s) )>0$, for $x\in\R$
and for each $s\in S$,
%
%e2.19 #&#
%
\begin{equation}
\label{GPcartesiandf>c} P\bigl(W(s)> x|W(s)>\omega_0\bigr)=\cases{
1, &\quad $x \leq\omega_0$,
\cr
\omega_0/x, &\quad $x >
\omega_0$.}
\end{equation}
\end{enumerate}
\end{cor}

\begin{pf} For (\ref{GPcartesiandf>}), similarly to the proof of
Proposition \ref{Wdf1prop},
\begin{eqnarray*}
P(W> w)&=&P \biggl(Y\geq\sup_{s\in S} \frac{w(s)}{V(s)}
\mbox{ and } \inf_{s\in S} V(s)>0 \biggr)
\\
&=&\int_{B_1} \inf_{s\in S}
\frac{v(s)}{w(s)}\mrmd \rho(v)=\rho(B_1)E \biggl(\inf
_{s\in S}\frac{V(s)}{w(s)} \Big| V\in B_1 \biggr).
\end{eqnarray*}

For (\ref{P>w>0}),
\[
P(W> w)=P(YV>w)=P \biggl(Y>\sup_{s\in S} \frac{w(s)}{V(s)}
\biggr)=E \biggl(\inf_{s\in S}\frac{V(s)}{w(s)} \biggr),
\]
using $Y$ standard Pareto and independent of $V$.

For (c) note that
\begin{eqnarray*}
P (W> w_0 )&=&P \Bigl(Y\inf_{s\in S} V(s)>
\omega_0 \Bigr)=E\min\biggl(1,\frac{\inf_{s\in S}V(s)}{\omega_0} \biggr)
\\
&=&\frac1{\omega_0}E\inf_{s\in S}V(s)>0.
\end{eqnarray*}
Then (\ref{GPcartesiandf>b}) follows from (\ref{P>w>0}).

For (d) note that, if $x>\omega_0$,
%
%e2.20 #&#
%
\begin{eqnarray}
\label{marginalsurvf} P \bigl(W(s)> x \bigr)&=&P \bigl(Y V(s)>x \bigr
)=E\min
\biggl(1,\frac
{V(s)}x \biggr)\nonumber\\[-9pt]\\[-8pt]
&=&x^{-1}EV(s)>0.\nonumber\\[-12pt]\nonumber
\end{eqnarray}
\upqed\end{pf}

Relation (\ref{GPcartesiandf>c}) indicates that one-dimensional
marginals, conditional on the process being larger than $\omega_0$,
behave like Pareto; a similar observation has been done by Rootz{\'e}n
and Tajvidi \cite{RootzenTajvidi} in the context of lower-dimensional
distributions.

Let $s_1,s_2\in S$ and $x>\omega_0$. From (\ref{marginalsurvf}),
\[
P \bigl(W(s_i)>x \bigr)=\frac{E (V(s_i) )}x >0,\qquad i=1,2,
\]
and, similarly
\[
P \bigl( W(s_1)>x,W(s_2)>x \bigr)=\frac{E (V(s_1)\wedge
V(s_2) )}x.
\]
Hence the statement $P ( W(s_1)>c,W(s_2)>c )=P
(W(s_1)>c )P (W(s_2)>c )$ for all $c>\omega_0$ is
equivalent to the statement $E (V(s_1)\wedge V(s_2)
)=c^{-1}E (V(s_1) )E (V(s_2) )$ for all
$c>\omega_0$, which is impossible. That is, independence in the Pareto
process between any two points is impossible.

For later use, we define next max-stable processes and give a
well-known property.

%
%de2.2 #&#
%
\begin{defi}\label{maxstabdef}
A process $\eta=\{\eta(s)\}_{s\in\R}\in C(\R)$ with non-degenerate
mar\-gi\-nals is called max-stable if, for $\eta_1, \eta_2, \ldots\,$,
i.i.d. copies of $\eta$, there are real continuous functions $c_n=\{
c_n(s)\}_{s\in\R}>0$ and $d_n=\{d_n(s)\}_{s\in\R}$ such that,
\[
\max_{1\leq i\leq n}\frac{\eta_i - d_n}{c_n}\stackrel{d} {=}\eta\qquad\mbox{for
all } n=1,2, \ldots.
\]
The process is called simple if its marginal distributions are standard
Fr\'echet, and then it will be denoted by $\bar\eta$.
\end{defi}

%
%pr2.3 #&#
%
\begin{prop}[(Penrose \cite{Penrose92}, Theorem 5)]\label{maxstabrepres}
All simple max-stable processes can be generated in the following way.
Consider a Poisson point process on $(0, \infty]$ with mean measure
$r^{-2}$ $dr$. Let
$ \{Z_i \}_{i=1}^{\infty}$ be a realization of this point
process. Further
consider i.i.d. stochastic processes $V_1, V_2,\ldots$ in $C^+(\R)$
with $E V_1(s)=1$
for all $s \in\R$ and $E \sup_{s\in\R} V(s)<\infty$. Then
%
%e2.21 #&#
%
\begin{equation}
\label{maxstabspectralV} \bar\eta=^d\max
_{i=1,2,\ldots} Z_i V_i.
\end{equation}
Conversely, each process with this representation is simple max-stable
(and one can take $V$ such that $\sup_{s\in\R}V(s)=c$ a.s. with $c>0$).
\end{prop}

Note that $\bar\eta$ depends on infinitely many processes $V_i$
whereas $W$ depends on just one of those processes (Theorem \ref
{Wtheor}(3)).

%s2.2 #&#
\subsection{The finite-dimensional setting}
\label{MSPpsect}

The theory of simple Pareto random \emph{vectors} (r.v.) can be
obtained from the results of Theorem~\ref{Wtheor}, by taking a
discrete set for $S$, $S=\{s_1,\ldots,s_d\}$ say. Consequently,
consider for this section the r.v. $(W_1,\ldots,W_d)=(W(s_1),\ldots,W(s_d))$.

%
%de2.3 #&#
%
\begin{defi}\label{Wddefi}
The r.v. $(W_1,\ldots,W_d)\in\R^d_+$ with threshold parameter
$\omega_0$ is called simple Pareto. The probability measure
%
%e2.22 #&#
%
\begin{equation}
\label{Qfindimdef} \rho(B)=P \biggl(\frac{\omega_0
(W_1,\ldots,W_d)}{\max_{i=1,\ldots,d} W_i}\in B \biggr)
\end{equation}
for $B\in{\cal B} (\{(w_1,\ldots,w_d)\in\R^d_+\dvt\max
(w_1,\ldots,w_d)=\omega_0\} )$ is again called the spectral measure.
\end{defi}

It follows again that for having all marginals Pareto, one would need
$\max(V(s_1),\ldots,V(s_n))=\omega_0$, for all $s_1,\ldots,s_n\in
S$ and all $n=1,\ldots,d$, which corresponds to $V\equiv\omega_0$
a.s., that is, the complete dependence case.

Nonetheless, we see that it is possible that some finite-dimensional
marginals of a Pareto process have a Pareto distribution. For example,
consider a situation where the maximum of the process occurs a.s. at
some fixed locations in $S$, $s_1,\ldots,s_d$ say. Then
$(W(s_1),\ldots,W(s_d))$ is a $d$-dimensional simple Pareto r.v. with
threshold parameter $\omega_0$. Moreover, any
$(W(s'_1),\ldots,W(s'_D))$ where
$\{s_1,\ldots,s_d\}\subset\{s'_1,\ldots,s'_D\}$ is a $D$-dimensional
simple Pareto r.v. with the same threshold parameter $\omega_0$.

One can give formulas for distribution functions, following similar
reasoning as before. The statement corresponding to Proposition \ref
{Wdf1prop} is
%
%e2.23 #&#
%
\begin{equation}
\label{lem21findim} P(W_1\leq w_1,\ldots,W_d\leq w_d)= \cases{\displaystyle  \rho
(B_0) \biggl
\{1-E \biggl(\max_{i\in\bar I_0}\frac{V_i}{w_i} \Big| B_0
\biggr) \biggr\}, &\quad if $\rho(B_0)>0$,
\cr
0, &\quad if $
\rho(B_0)=0$,}
\end{equation}
where $I_0=\{1\leq i\leq d\dvt w_i=0\}$, $\bar I_0=\{1\leq i\leq d\dvt
w_i\neq0\}$ and
$B_0=\{(V_1,\ldots,V_d)\dvt V_i=0$ for $i\in I_0$ and
$\min_{i\in\bar I_0}\frac{w_i}{V_i}\geq1\}$.

The statement corresponding to Proposition \ref{Wdfprop}, with
$w_i>0$ for all $i=1,\ldots,d$, is
%
%e2.24 #&#
%
\begin{equation}
\label{lem22findim} P(W_1\leq w_1,\ldots,W_d\leq w_d)=E \biggl(\max
_{1\leq i\leq d}
\frac
{V_i}{w_i\wedge\omega_0} \biggr)-E \biggl(\max_{1\leq i\leq d}\frac
{V_i}{w_i}
\biggr),
\end{equation}
which corresponds to formula (2) in Definition 2.1 from Rootz{\'e}n and
Tajvidi \cite{RootzenTajvidi}.

Note that Rootz\'en and Tajvidi's formula, (2) in Definition 2.1, only
holds for a vector $(x,y)$ -- for simplicity, we take bivariate vectors
in the following discussion -- larger than the vector of the lower
endpoints of the marginal distributions. The following example
illustrates this fact.

Apply Rootz\'en and Tajvidi's formula with
$G(x,y)=\RMe^{-((x+1)^{-1}+(y+1)^{-1})}$, for $x>-1$ and $y>-1$, that is,
the r.v. constructed from two independent unit Fr\'echet random
variables shifted by $-1$. Then, their formula for the Pareto r.v.
shifted by $(1,1)$ is,
%
%e2.25 #&#
%
\begin{equation}
\label{exfindimRT} H(x,y)= \cases{\displaystyle  \frac1 2 \biggl(\frac1{x\wedge
0}-\frac1{y
\wedge0}-\frac1{x}-\frac1{y} \biggr)\vspace*{2pt}\cr
\hspace*{19.65pt}\mbox{\quad $(x>0 \mbox{ and } y\geq1)
\mbox{ or }
(y>0 \mbox{ and } x\geq1)$},
\cr
0, \qquad 0<x\leq1, 0<y\leq1.}
\end{equation}
Now note that this does not properly accommodate the positive mass that
exists on the axis.

Our alternative approach leads, in this case, to the following
distribution function. Consider the bivariate Pareto r.v. $(YB,
Y(1-B))$, with $B$ Bernoulli $(1/2)$. Then, by direct calculations or by
applying (\ref{lem21findim}) one obtains the distribution function
%
%e2.26 #&#
%
\begin{equation}
\label{exfindim} P\bigl(YB\leq x,Y(1-B)\leq y\bigr)= \cases{\displaystyle  \frac1 2
\biggl( 2-
\frac1{x}-\frac1{y} \biggr), &\quad if $x\geq1,y\geq1$,
\vspace*{2pt}\cr
\displaystyle \frac1 2 \biggl( 1-
\frac1{x} \biggr), &\quad if $x\geq1,0\leq y< 1$,
\vspace*{2pt}\cr
\displaystyle \frac1 2 \biggl( 1-\frac1{y}
\biggr), &\quad if $y\geq1,0\leq x< 1$,
\vspace*{2pt}\cr
0, &\quad otherwise.}
\end{equation}
Regard that (\ref{exfindimRT}) and (\ref{exfindim}) are the same
except when $x=0$ or $y=0$.

Another remark on Rootz{\'e}n and Tajvidi \cite{RootzenTajvidi}: their
Theorem 2.2(ii) is not completely correct. It is not sufficient to
require condition (6) of the same paper for $x,y>0$. A counter example
is given by
\[
P(X>x \mbox{ or } Y>y)= \bigl(\tfrac1 2 \RMe^{-2(x\vee0)}+\tfrac1 2
\RMe^{-2(y\vee0)} \bigr)^{1/2},\qquad x\vee y \geq0,
\]
and zero elsewhere. This distribution satisfies (6) for $x,y>0$ but not
for all $(x,y)$ and it is not a generalized Pareto distribution.

%s2.3 #&#
\subsection{The generalized Pareto process}
\label{GSPpsect}
The more general processes with continuous extreme value index function
$\gamma=\{\gamma(s)\}_{s\in S}$, location and scale functions $\mu=\{
\mu(s)\}_{s\in S}$ and $\sigma=\{\sigma(s)\}_{s\in S}$ is defined as follows.

%
%de2.4 #&#
%
\begin{defi}\label{GPdefi2}
Let $W$ be a simple Pareto process, $\mu,\sigma,\gamma\in C(S)$ with
$\sigma>0$. The generalized Pareto process $W_{\mu,\sigma,\gamma
}\in C(S)$ is defined by,
%
%e2.27 #&#
%
\begin{equation}
W_{\mu,\sigma,\gamma}=\mu+\sigma\frac{W^{\gamma}-1}\gamma
\end{equation}
with all operations taken componentwise (recall the convention
explained in the end of Section~\ref{Introdsect}).
\end{defi}

The result corresponding to Corollary \ref{supxcor} is the following.

%
%co2.5 #&#
%
\begin{cor}\label{stparetoWgcor}
The random variable $\sup_{s\in S} \{ (1+\gamma(s)\frac
{W_{\mu,\sigma,\gamma}(s)-\mu(s)}{\sigma(s)} )^{1/\gamma
(s)} \}\omega_0^{-1}$ has standard Pareto distribution.
\end{cor}

%change to $P(W\in rA)=r^{-\alpha}P(W\in A)$, for all $r>1$ and $A\in{
%W(s)/\omega_0$ would have the same distribution as $Y$.\end{rem}

The process satisfies the following stability property.

%
%pr2.4 #&#
%
\begin{prop}\label{stparetoWgprop}
For any generalized Pareto process $W_{\mu,\sigma,\gamma}$,
%
%e2.28 #&#
%
\begin{equation}
\label{homogWg} P \biggl( \biggl(1+\gamma\frac{W_{\mu,\sigma,\gamma}-\mu
}{\sigma
}
\biggr)^{1/\gamma}\in rA \biggr) =r^{-1}P \biggl( \biggl(1+\gamma
\frac{W_{\mu,\sigma,\gamma}-\mu
}{\sigma} \biggr)^{1/\gamma}\in A \biggr)
\end{equation}
for all $r>1$ and $A\in{\cal B} (C_{\omega_0}^{+}(S) )$.
Moreover, there exist normalizing functions $u(r)$ and $s(r)$ such that
%
%e2.29 #&#
%
\begin{eqnarray}
\label{stabWg} &&P \biggl( \biggl(1+\gamma\frac{W_{\mu,\sigma,\gamma
}-u(r)}{s(r)}
\biggr)^{1/\gamma}\in A \Big| \sup_{s\in S} \biggl(1+\gamma
\frac{W_{\mu,\sigma,\gamma
}-u(r)}{s(r)} \biggr)^{1/\gamma}>\omega_0 \biggr)
\nonumber\\[-8pt]\\[-8pt]
&&\quad=P \biggl( \biggl(1+\gamma\frac{W_{\mu,\sigma,\gamma}-\mu}{\sigma
} \biggr)^{1/\gamma}\in A \biggr)
\nonumber
\end{eqnarray}
for all $r>1$ and $A\in{\cal B} (C_{\omega_0}^{+}(S) )$.

Conversely, if (\ref{stabWg}) holds and $\sup_{s\in S} \{
(1+\gamma(s)\frac{W_{\mu,\sigma,\gamma}(s)-\mu(s)}{\sigma
(s)} )^{1/\gamma(s)} \}\omega_0^{-1}$ has a standard
Pareto distribution, then (\ref{homogWg}) holds.
\end{prop}

\begin{pf}
Relation (\ref{homogWg}) is direct from Definition \ref{GPdefi2} and
(\ref{GPhomogeneity}). Then, with $u(r)=\mu+\sigma(r^\gamma
-1)/\gamma$ and $s(r)=\sigma r^\gamma$, relation (\ref{stabWg}) is
easily shown to be true by (\ref{homogWg}) and Corollary \ref
{stparetoWgcor}.

Conversely, for all $r>1$ and $A\in{\cal B} (C_{\omega
_0}^{+}(S) )$,
\begin{eqnarray*}
&&\frac{P ( (1+\gamma({W_{\mu,\sigma,\gamma
}-u(r)})/{s(r)} )^{1/\gamma}\in A )}{P (\sup_{s\in
S} (1+\gamma({W_{\mu,\sigma,\gamma}-u(r)})/{s(r)}
)^{1/\gamma}\omega_0^{-1}>1 )}
\\
&&\quad=\frac{P ( (1+\gamma({W_{\mu,\sigma,\gamma}-\mu
})/{\sigma} )^{1/\gamma}\in r A )}{r^{-1}} \\
&&\quad=P \biggl( \biggl(1+\gamma\frac
{W_{\mu,\sigma,\gamma}-\mu}{\sigma}
\biggr)^{1/\gamma}\in A \biggr)
\end{eqnarray*}
by (\ref{stabWg}) and $\sup_{s\in S} \{ (1+\gamma(s)\frac
{W_{\mu,\sigma,\gamma}(s)-\mu(s)}{\sigma(s)} )^{1/\gamma
(s)} \}\omega_0^{-1}$ being standard Pareto distributed.
\end{pf}

The result corresponding to Proposition \ref{Wdfprop} on distribution
functions is now, for $w>0$:
\begin{eqnarray*}
P(W_{\mu,\sigma,\gamma}\leq w)&=&E \biggl\{\sup_{s\in S}V(s) \biggl(
\biggl(1+\gamma(s)\frac{w(s)-\mu(s)}{\sigma(s)} \biggr)^{1/\gamma
(s)}\wedge
\omega_0 \biggr)^{-1} \biggr\}
\\
&&{}- E \biggl\{\sup_{s\in S}V(s) \biggl(1+\gamma(s)
\frac{w(s)-\mu
(s)}{\sigma(s)} \biggr)^{-1/\gamma(s)} \biggr\}
\end{eqnarray*}
for $1+\gamma(w-\mu)/\sigma\in C^+(S)$.
%P(W_{\mu,\sigma,\gamma}>w)=E\left\{\inf_{s\in\R}V(s)\left(1+\gamma(s)
% and for the univariate marginals we have,
%P(W(s)> x(s))=E\left\{V(s)\left(1+\gamma(s)(x(s)-\mu(s))/\sigma(s)

%s3 #&#
\section{Domain of attraction}
\label{domattrsect}

Let us start with the characterization of the domain of attraction of a
max-stable process. This result will lead directly to a
characterization of the domain of attraction of a generalized Pareto
process. The following is a slight variation and extension of Theorem
9.5.1 of de Haan and Ferreira \cite{HaanFerreira}.

Denote by $\bar\eta=\{\bar\eta(s)\}_{s\in S}$ any simple max-stable
process in $C^+(S)$ (cf. Definition \ref{maxstabdef}). Any max-stable
process $\eta=\{\eta(s)\}_{s\in S}$ in $C(S)$ can be represented by
$\eta=(\bar\eta^\gamma-1)/\gamma$, for some $\bar\eta$ and
continuous function $\gamma=\{\gamma(s)\}_{s\in S}$. For simplicity,
we always take here
\[
C_1^+(S)=\Bigl\{f \in C^+(S)\dvt\sup_{s\in S} f(s)\geq1
\Bigr\}
\]
that is, w.l.g. consider the constant $\omega_0$ introduced in
Section \ref{SPpsect} equal to 1.
For $X$ a random element of $C(S)$, suppose the marginal distribution
functions $F_s(x)=P(X(s)\leq x)$ are continuous in $x$, for all $s\in S$.

%
%th3.1 #&#
%
\begin{theo}\label{maxdomtheo}
Suppose $X,X_1, X_2, \ldots$ are i.i.d. random elements of $C(S)$. The
following statements are equivalent.
\begin{enumerate}
\item There exists a max-stable stochastic process $\eta\in C(S)$ with
continuous index function $\gamma$, and $a_n>0$ and $b_n$ in $C(S)$
such that
%
%e3.1 #&#
%
\begin{equation}
\label{max-domattr} \biggl\{\max_{1\leq i\leq n} \frac{X_i(s) - b_n(s)}{a_n(s)}
\biggr\} _{s \in S}\to^d \bigl\{\eta(s) \bigr\}_{s \in S}
\end{equation}
in $C(S)$ ($\to^d$ denotes weak convergence or convergence in
distribution). The normalizing functions are w.l.g. chosen in such a
way that $-\log P(\eta(s)\leq x)=(1+\gamma(s)x)^{-1/\gamma(s)}$ for
all $x$ with $1+\gamma(s)x>0$, $s\in S$.
\item There exist continuous functions $\gamma$, $a_t>0$ and $b_t$
such that
%
%e3.2 #&#
%
\begin{equation}
\label{basicE29} \lim_{t\to\infty} tP \biggl(\frac{X(s)- b_t (s)}{a_t (s)}>x
\biggr)=\bigl(1+\gamma(s)x\bigr)^{-1/\gamma(s)},\qquad 1+\gamma(s)x>0,
\end{equation}
uniformly for $s\in S$ and, for the normalized process
\[
T_tX= \biggl(1+\gamma\frac{X-b_t}{a_t} \biggr)_+^{1/\gamma}
\]
we have
%
%e3.3 #&#
%
\begin{equation}
\label{WDomAttreqb} \lim_{t\to\infty} \frac{P (\sup_{s\in S}T_tX(s)>x
)}{P (\sup_{s\in S}T_tX(s)>1 )}=\frac1 x
\qquad\mbox{for all } x>1
\end{equation}
and
%
%e3.4 #&#
%
\begin{equation}
\label{WDomAttreqc} \lim_{t\to\infty} P \biggl(\frac{T_tX}{\sup_{s\in
S}T_tX(s)}\in B
\Big|\sup_{s\in S}T_tX(s)>1 \biggr)=\rho(B)
\end{equation}
for each $B\in{\cal B} (\bar C_1^+(S) )$ with $\rho
(\partial B)=0$, with $\rho$ some probability measure on $\bar C_1^+(S)$.
\end{enumerate}
\end{theo}

The following shows that the same conditions are valid for the domain
of attraction of a generalized Pareto process.

%
%th3.2 #&#
%
\begin{theo}\label{RTdirectstatcor}
1. The conditions of Theorem \ref{maxdomtheo} imply
\[
\lim_{t\to\infty} P \Bigl(T_tX\in A \Big|\sup
_{s\in
S}T_tX(s)>1 \Bigr)=P(W\in A)
\]
with $A\in{\cal B} (C_1^+(S) )$, $P(\partial A)=0$ and W
some simple Pareto process.

2. Conversely suppose that there exists a function $\tilde b_u=\{
\tilde b_u(s)\}_{s\in S}$, that is continuous in $s$ for each $u$ and
increasing in $u$, and with the property that $P(X(s)>\tilde
b_u(s)$ for some $s\in S)\to0$ as $u\to\infty$, and a
continuous function (in $s$), $\tilde a_u=\{\tilde a_u(s)\}_{s\in S}>0$
such that, for some probability measure $\tilde P$ on ${\cal B} (
C(S) )$,
\[
\lim_{u\to\infty}P \biggl(\frac{X-\tilde b_u}{\tilde a_u}\in A \Big|
X(s)-\tilde
b_u(s)>0 \mbox{ for some } s\in S \biggr) =\tilde P(A)
\]
for all $A\in{\cal B} (C(S) )$ and $\tilde P(\partial
A)=0$. Then the results of Theorem \ref{maxdomtheo} hold.
\end{theo}

\begin{pf}
Statement 1 follows directly from Theorem \ref{maxdomtheo}.

We prove statement 2:

By the conditions on $\tilde b_u$, we can determine $q=q(t)$ such that
$P(X(s)>\tilde b_{q(t)}(s)\mbox{ for some } s\in S)=1/t$. Then with
$b_t(s)=\tilde b_{q(t)}(s)$ and $a_t(s)=\tilde a_{q(t)}(s)$,
\[
\lim_{t\to\infty}tP \biggl(\frac{X-b_t}{a_t}\in C \mbox{ and } X(s)>
b_t(s) \mbox{ for some } s\in S \biggr)=\tilde P(C)
\]
for all $C\in{\cal B} (C(S) )$ and $\tilde P(\partial
C)=0$. In particular, if $\inf\{\sup_{s\in S}f(s)\dvt f\in C\}>0$ we have
%
%e3.5 #&#
%
\begin{equation}
\label{limtP} \lim_{t\to\infty}tP \biggl(\frac{X-b_t}{a_t}\in C
\biggr)=\tilde P(C).
\end{equation}

We proceed as usual in extreme value theory. Fix for the moment $s\in
S$. It follows that for $x>0$
\[
\lim_{t\to\infty}tP \bigl(X(s)>b_t(s)+xa_t(s)
\bigr)=\tilde P\bigl\{ f\dvt f(s)>x\bigr\}.
\]
%
%Note that this is not convergence in distribution. Nevertheless
Let $U_s$ be the inverse function of $1/P(X(s)>x)$ and $V(s)$ be the
inverse function of $1/\tilde P \{f\dvt f(s)>x \}$. Then
\[
\lim_{t\to\infty}\frac{U_{tx}(s)-b_t(s)}{a_t(s)}=V_x(s)\qquad \mbox{for }
x>0.
\]
It follows (Lemma 10.4.2, p. 340, in de Haan and Ferreira \cite
{HaanFerreira}) that for some real $\gamma(s)$ and all $x>0$
%
%e3.6 #&#
%
\begin{equation}
\label{limtPmarg} \lim_{t\to\infty}\frac{b_{tx}(s)-b_t(s)}{a_t(s)}=
\frac{x^{\gamma
(s)}-1}{\gamma(s)} \quad\mbox{and}\quad \lim_{t\to\infty}\frac
{a_{tx}(s)}{a_t(s)}=x^{\gamma(s)}.
\end{equation}
Since the limit process has continuous paths, the function $\gamma$
must be continuous on $S$.

Now replace $t$ in (\ref{limtP}) by $ct$ where $c>0$. Then
\[
\lim_{t\to\infty}tP \biggl(\frac{b_t(s)-b_{tc}(s)}{a_{tc}(s)}+\frac
{a_t(s)}{a_{tc}(s)}
\frac{X-b_t}{a_t}\in C \biggr)=\frac1 c\tilde P(C)
\]
hence, by (\ref{limtPmarg})
\[
\lim_{t\to\infty}tP \biggl( \biggl(1+\gamma\frac{X-b_t}{a_t}
\biggr)^{1/\gamma}\in c (1+\gamma C )^{1/\gamma} \biggr)=\frac1 c\tilde P(C)
\]
and by (\ref{limtP})
\[
\lim_{t\to\infty}tP \biggl( \biggl(1+\gamma\frac{X-b_t}{a_t}
\biggr)^{1/\gamma}\in(1+\gamma C )^{1/\gamma} \biggr)=\tilde P(C).
\]
Write $P(A)=\tilde P ( (A^{\gamma}-1 )/\gamma
)$. Then
\[
\lim_{t\to\infty}tP (T_tX\in A )=P(A)
\]
with $P(cA)=c^{-1}P(A)$, for all $c>0$ and $A\in{\cal B}
(C(S) )$ such that $\inf\{\sup_{s\in S}f(s)\dvt f\in A\}>1$ and
$P(\partial A)=0$.
The rest is like the proof of the equivalence between (2b) and (2c) of
Theorem~9.5.1 in de Haan and Ferreira \cite{HaanFerreira}.
\end{pf}

%
%ex3.1 #&#
%
\begin{ex}
Any max-stable process is in the domain of attraction of a generalized
Pareto process, with $\rho$ given by the probability measure of $V$ from
(\ref{maxstabspectralV}).
\end{ex}

%
%ex3.2 #&#
%
\begin{ex}
Any Pareto process with spectral measure $\rho$ is in the domain of
attraction of a max-stable process where the underlying process $V$
(cf. representation (\ref{maxstabspectralV})) has probability measure
$\rho$.
\end{ex}

%
%ex3.3 #&#
%
\begin{ex}
The finite-dimensional distributions of the moving maximum processes
obtained in de Haan and Pereira \cite{HaanPereira06} can be applied to
obtain the finite-dimensional distributions of the corresponding Pareto process.
\end{ex}

%
%ex3.4 #&#
%
\begin{ex}[(Regular variation (de Haan and Lin \cite{HaanLin}, Hult and
Lindskog \cite{HultLindskog}))]
A stochastic process $X$ in $C(S)$ is regularly varying if and only if
there exists an $\alpha>0$ and a probability measure $\rho$ such that,
%
%e3.7 #&#
%
\begin{equation}
\label{RV} \frac{P (\sup_{s\in S}X(s)>tx, X/\sup_{s\in S}X(s)\in\cdot
)}{P (\sup_{s\in S}X(s)>t )}\to^d x^{-\alpha} \rho(\cdot),\qquad x>0,
t\to\infty,
\end{equation}
on $\{f\in C(S)\dvt\sup_{s\in S}f(s)=1\}$. Hence, a regularly varying
process such that (\ref{basicE29}) holds for the marginals, satisfies
the conditions of Theorem \ref{maxdomtheo}, in particular with
$\gamma=1/\alpha$, $b_t=t$ and $a_t=t/\alpha$; note that the index
function is constant.

On the other hand, the normalized process $T_tX$, with $T_tX$
satisfying (\ref{WDomAttreqb})--(\ref{WDomAttreqc}), is regularly
varying with $\alpha=1$ and spectral measure $\rho$ on $\bar C_1^+(S)$.
\end{ex}

%
%re3.1 #&#
%
\begin{rem}
As seen in Section \ref{MSPpsect}, our analysis is also valid in the
finite-dimensional set-up. The main difference from Rootz{\'e}n and
Tajvidi \cite{RootzenTajvidi} is that their analysis is entirely based
on distribution functions whereas ours is more structural. Here are
some remarks on their domain of attraction results.

Let $\bar F=1-F$ with $F$ some $d$-variate distribution function,
$\mathbf x=(x_1,x_2,\ldots,x_d)\in\R^d$, and $\mathbf u(\cdot
)=(u_1(\cdot),u_2(\cdot),\ldots,u_d(\cdot))$ and $\bolds{\sigma
}$ the normalizing functions considered in Rootz{\'e}n and Tajvidi
\cite{RootzenTajvidi} (see, e.g., their definition of $\mathbf
{X_u}$). By using $\bolds\sigma(\mathbf x t)/\bolds\sigma(t)\to
(x_1^{\gamma_1},x_2^{\gamma_2},\ldots,x_d^{\gamma_d} )$
and
$ (\mathbf u(\mathbf x t)-\mathbf u(t) )/\bolds\sigma
(t)\to(\frac{x_1^{\gamma_1}-1}{\gamma_1},\frac{x_2^{\gamma
_2}-1}{\gamma_2},\ldots,\frac{x_d^{\gamma_d}-1}{\gamma_d}
)$, $t\to\infty$,\vspace*{2pt} for some reals $\gamma_1,\gamma_2,\ldots,\gamma
_d$ (cf. proof of Theorem 2.1(ii) in Rootz{\'e}n and Tajvidi \cite
{RootzenTajvidi}) and by
\[
\bar F^*(\mathbf x):=\bar F \bigl(u_1(x_1),u_2(x_2),\ldots,u_d(x_d)
\bigr),
\]
one simplifies their relation (19) to
\[
t\bar F^*(t\mathbf x)\to-\log G \biggl(\frac{x_1^{\gamma_1}-1}{\gamma
_1},\frac{x_2^{\gamma_2}-1}{\gamma_2},\ldots,\frac{x_d^{\gamma
_d}-1}{\gamma_d} \biggr),
\]
and one simplifies their relation (6) to
\[
P \bigl(\mathbf X^*\leq t\mathbf x|\mathbf X^*\nleq t\mathbf1 \bigr)=P
\bigl(
\mathbf X^*\leq\mathbf x \bigr)
\]
for $t\geq1$. Hence, one can take
$\mathbf u(t):= (\frac{t^{\gamma_1}-1}{\gamma_1},\frac
{t^{\gamma_2}-1}{\gamma_2},\ldots,\frac{t^{\gamma_d}-1}{\gamma
_d} )$ in Theorem 2.2 of that paper.
\end{rem}

%s4 #&#
\section{View towards application and simulation}
\label{applicsect}

%s4.1 #&#
\subsection{Towards application}
Suppose the domain of attraction condition (\ref{nudefintrod}) holds.
Define $B=\{f\in C^+(S)\dvt\sup_{s\in S} f(s)>1\}$. Let $A$ be a Borel
set in $C^+(S)$. Then applying (\ref{nudefintrod}) twice we get
\[
\lim_{t\to\infty} P (T_tX\in A |T_tX\in B
)=\frac
{\nu(A\cap B)}{\nu(B)}=P(W\in A)
\]
with $W$ a simple Pareto process. This is the content of Theorem \ref
{RTdirectstatcor}(1) and it gives the basis to the
peaks-over-threshold method in function space, as it gives a limit
probability distribution on~$B$.

A similar reasoning holds with $B$ replaced by a different set $B'$ as
long as $\inf\{\sup_{s\in S} f(s)\dvt f\in B'\}>0$. Consider in
particular $B'=\{f\in C^+(S)\dvt\max_{i=1,\ldots, p} T_t f(s_i)\geq1\}
$, for some integer $p$. Then
%
%e4.1 #&#
%
\begin{equation}
\label{limnuappl} \lim_{t\to\infty} P \bigl(T_tX\in A
|T_tX\in B' \bigr)=\frac
{\nu(A\cap B')}{\nu(B')}=
\frac{P(W\in A\cap B')}{P(W\in B')},
\end{equation}
which is a generalized Pareto distribution.

Now, we proceed as in the peaks-over-threshold method for scalar
observations: let $k=k(n)$ be a sequence of integers with $\lim_{n\to
\infty} k(n)=\infty$ and $\lim_{n\to\infty} k(n)/n=0$, as $n\to
\infty$. Suppose that we have $n$ independent observations of the
process $X$ in the domain of attraction. Select those observations
satisfying $X(s_i)>b_{n/k}(s_i)$, for some $i=1,2,\ldots,p$. The
probability distribution of those selected observations is
approximately the right-hand side of (\ref{limnuappl}), that is,
generalized Pareto. This seems a useful applicable form of the
peaks-over-threshold method in this framework as it suggests estimating
the spectral measure using observations that exceed a threshold at some
discrete points in the space only.

%s4.2 #&#
\subsection{Towards simulation}\label{sec4.2}
`Deltares' is an advisory organization of the Dutch government
concerning (among others) the safety of the coastal defenses against
severe wind storms. One studies the impact of severe storms on the
coast, storms that are so severe that they have never been observed. In
order to see how these storms look like, it is planned to simulate wind
fields on and around the North Sea using certain climate models. These
climate models simulate independent and identically distributed
(i.i.d.) wind fields similar to the ones that could be observed (but
that are only partially observed). Since the model runs during a
limited time, some of the wind fields will be connected with storms of
a certain severity but we do not expect to see really disastrous storms
that could endanger the coastal defenses. The question put forward by
Deltares is: can we get an idea how the really disastrous wind fields
look like on the basis of the `observed' wind fields? We want to show
that this can be done using the generalized Pareto process.

Consider a continuous stochastic processes $\{X(s)\}_{s\in S}$ where
$S$ is a compact subset of $\R^d$. Suppose that the probability
distribution of the process is in the domain of attraction of some
max-stable process, that is, there exist functions $a_n>0$ and $b_n$
such that the sequence of i.i.d. processes
\[
\biggl\{\max_{1\leq i\leq n} \frac{X_i(s) - b_n(s)}{a_n(s)} \biggr\}
_{s \in S}
\]
converges to a continuous process, say $\eta$, in distribution in
$C(S)$. Then $\eta$ is a max-stable process.

Define
\begin{eqnarray*}
T_tX(s)&:=& \biggl(1+\gamma(s)\frac{X(s)-b_t(s)}{a_t(s)}
\biggr)_+^{1/\gamma(s)},
\\
R_{T_tX}&:=&\sup_{s\in S} T_tX(s).
\end{eqnarray*}
Then
\[
T_t^{\leftarrow}f(s):=a_t(s)\frac{(f(s))^{\gamma(s)}-1}{\gamma
(s)}+b_t(s)
\qquad\mbox{for } f\in C^+(S).
\]
As before, suppress the $s$ from now on. Then, with $t_0$ some large constant,
%
%e4.2 #&#
%
\begin{eqnarray}\label{condprobTT}
&&P \biggl(\frac{T_t^{\leftarrow}t_0 T_t
X-b_{tt_0}}{a_{tt_0}}\in A \Big|R_{T_tX}>1 \biggr)
\nonumber\\
&&\quad=P \biggl\{\frac{a_t}{a_{tt_0}}\frac{ [t_0 (1+\gamma
({X-b_t})/{a_t} )^{1/\gamma} ]^{\gamma}-1}{\gamma
}-\frac{b_{tt_0}-b_t}{a_{tt_0}}\in A
\Big|R_{T_tX}>1 \biggr\}
\nonumber
\\
&&\quad=P \biggl\{\frac{a_t t_0^{\gamma}}{a_{tt_0}}\frac{1+\gamma
({X-b_t})/{a_t}-t_0^{-\gamma}}{\gamma}-\frac{b_{tt_0}-b_t}{a_{tt_0}}\in A
\Big|R_{T_tX}>1 \biggr\}
\\
&&\quad=P \biggl\{\frac{a_t t_0^{\gamma}}{a_{tt_0}} \biggl(\frac
{X-b_t}{a_t}-t_0^{-\gamma}
\biggl[\frac{b_{tt_0}-b_t}{a_{t}}-\frac
{t_0^{\gamma}-1}{\gamma} \biggr] \biggr)\in A
\Big|R_{T_tX}>1 \biggr\}
\nonumber
\\
&&\quad=P \biggl\{\frac{X-b_t}{a_t}\in\frac{a_{tt_0} t_0^{-\gamma}}{a_t}
A+t_0^{-\gamma}
\biggl(\frac{b_{tt_0}-b_t}{a_{t}}-\frac{t_0^{\gamma
}-1}{\gamma} \biggr)  \Big|R_{T_tX}>1
\biggr\}.
\nonumber
\end{eqnarray}
Since,
\[
\frac{a_{tt_0}(s) t_0^{-\gamma(s)}}{a_t(s)}\to1 \quad\mbox{and}\quad \frac
{b_{tt_0}(s)-b_t(s)}{a_{t}(s)}-\frac{t_0^{\gamma
(s)}-1}{\gamma(s)}\to0
\]
uniformly for $s\in S$, the limit of this probability, as $t\to\infty
$, is the same as the limit of
%
%e4.3 #&#
%
\begin{equation}
\label{condprobX} P \biggl(\frac{X-b_t}{a_t}\in A \Big|R_{T_tX}>1 \biggr),
\end{equation}
which is generalized Pareto by Theorem \ref{RTdirectstatcor}(1).

In this subsection, we are not so much interested in estimating the
joint limit distribution (which is the peaks-over-threshold method) but
in the fact that the two conditional distributions (\ref{condprobTT})
and (\ref{condprobX}) are approximately the same.

Suppose, for example, that we have observed wind fields over a certain
area during some time. Then we are likely to find some rather heavy
storms, that is, ones that satisfy $X\nleq b_n$. These are the moderately
heavy storms. However, we want to know how the storm field of a really
heavy storm (i.e., $X\nleq b_N$ with $N>n$) looks like. That is exactly
what relation (\ref{condprobTT}) does. Take a moderately heavy storm
$X$ and transform it to $T_n^{\leftarrow}\frac N n T_n X$. This
results in the storm field of a really heavy storm by relation (\ref
{condprobTT}).

Notice then what we do here is similar to prediction or kriging, not
estimating a distribution function.

The reasoning above also holds with estimated functions of $\gamma$,
$a$ and $b$, on the basis of $k$th upper order statistics and taking $t=n/k$.

Under the above framework, we propose the following simulation method:
\begin{enumerate}[(5)]
\item[(1)] Let $X_1,X_2,\ldots,X_n$ be i.i.d. and let the underlying
distribution satisfy the conditions above, namely that the probability
distribution is in the domain of attraction of some max-stable process.
\item[(2)] Estimate the functions $\gamma$, $a$ and $b$ (de Haan and
Lin \cite{HaanLin2003}, Einmahl and Lin \cite{EinmahlLin}); denote
the estimators by $\hat\gamma$, $\hat a$ and $\hat b$. Note that this
procedure provides us with a number $k$ that reflects the threshold for
estimating the parameters.
\item[(3)] Select from the normalized processes
\[
\hat T_{n/k}X_i:= \biggl(1+\hat\gamma\frac{X_i-\hat b_{n/k}}{\hat
a_{n/k}}
\biggr)_+^{1/\hat\gamma},\qquad i=1,\ldots, n,
\]
those that satisfy $X_i(s)>\hat b_{n/k}(s)$ for some $s\in S$, that is,
for which $\hat R_{T_{n/k}X_i}>1$.
\item[(4)] Multiply these processes by a (large) factor $t_0$; this
brings the processes to a higher level without changing the
distribution essentially.
\item[(5)] Finally undo the normalization, that is, in the end we
obtain the processes
\[
\hat T_{n/k}^{\leftarrow}t_0 \hat T_{n/k}
X_i \mbox{ for those } X_i \mbox{ for which } \hat
R_{T_{n/k}X_i}>1.
\]
These processes are peaks-over-threshold processes with respect to a
much higher threshold (namely $b_{tt_0}$) than the processes $X_i$ for
which $\hat R_{T_{n/k}X_i}>1$ (with threshold $b_t$).
\end{enumerate}
%
%f1 #&#
%
\begin{figure}

\includegraphics[scale=0.97]{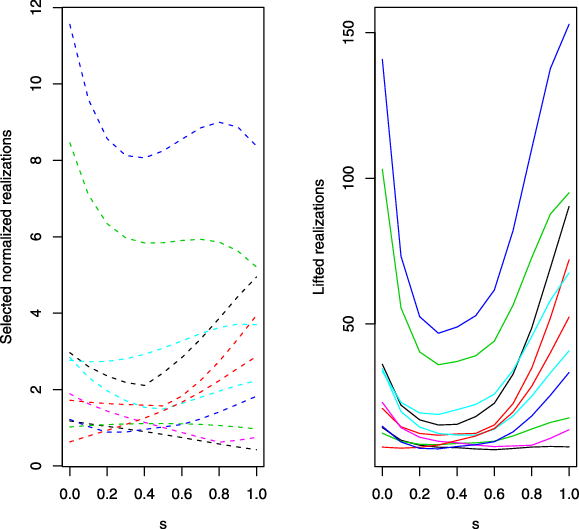}

\caption{(Left): the realizations $\hat T_{n/k}X_i$, for which
$\hat
R_{T_{n/k}X_i}>1$, obtained from the moving maximum process with
standard Gaussian density; (right) the lifted realizations $\hat
T_{n/k}^{\leftarrow}t_0 \hat T_{n/k} X_i$ with $t_0=10$.} \label{gauss1fig}
\end{figure}
%

%
%re4.1 #&#
%
\begin{rem}
Note that an alternative procedure under the maximum domain of
attraction condition would be, first to estimate the spectral measure
and then to simulate a generalized Pareto process from there. But the
estimation of the spectral measure is more difficult (de Haan and Lin
\cite{HaanLin}) although this procedure is less restrictive on the
number of observations that can be simulated.
\end{rem}

%s4.3 #&#
\subsection{Simulations}
We exemplify the lifting procedure with the process $X(s)=Z(s)^{\gamma
(s)}$, with $\gamma(s)=1-s(1-s)^2$, $s\in[0,1]$, and $Z$ is the
moving maximum process with standard Gaussian density. The $Z$ process
can be easily simulated in the R-package due to Ribatet \cite
{Ribatet}. Figure \ref{gauss1fig}\vadjust{\goodbreak} is represented by the 11 out of 20
realizations, normalized for which $\hat R_{T_{n/k}X_i}>1$, and lifted
ones $\hat T_{n/k}^{\leftarrow}t_0 \hat T_{n/k} X_i$ with $t_0=10$.

\section*{Acknowledgements}

An anonymous referee pointed out an error in the original version of
Proposition \ref{Wdfprop}; we are grateful for this.

The presented work in Section \ref{sec4.2} is part of the SBW (Strength and Load
on Water Defences) project, commissioned by Rijkswaterstaat Centre for
Water Management in the Netherlands and carried out by Deltares, The
Netherlands.

Research partially supported by FCT Project PTDC/MAT/112770/2009 and
FCT-PEst-OE/MAT/UI0006/2011.

%suskaldyti doi

% imsref loaded by lrinkeviciute, 2013-08-19 12:21:36
% imsref loaded by lrinkeviciute, 2013-08-19 12:26:12
%

\printhistory

\end{document}